\documentclass[twoside,11pt]{article}
 \usepackage{bbm}
 \usepackage{mathrsfs}
  \usepackage{amsfonts}
\usepackage{amsmath}
\newcommand{\proof}{\noindent {\bf Proof}}

\def\sq{\vee}

\begin{document}
  \title{The critical groups for $K_m\vee P_n$ and $P_m\vee P_n$
 \thanks{Supported by ``the Fundamental Research Funds for the Central Universities" and the NSF of the People's Republic of China(Grant
 No. 10871189).}
  }
 \author{Wei-Na Shi,\quad  Yong-Liang Pan\thanks{Corresponding author. Email: ylpan@ustc.edu.cn},\quad Jian Wang
 \\
  {\small Department of Mathematics, University of Science and Technology
         of China}\\
  {\small Hefei, Auhui 230026, The People's Republic of China}\\
}
\date{}
\maketitle {\centerline{\bf\sc Abstract}\vskip 8pt Let $G_1\sq G_2$
denote the graph obtained from $G_1+G_2$ by adding new edges from
each vertex of $G_1$ to every vertex of $G_2$. In this paper, the
critical groups of the graphs $K_m\sq P_n$$(n\geq4)$ and $P_m\sq
P_n$$(m\geq4,n\geq5)$ are determined.

\par \vskip 0.5pt {\bf Keywords}  Graph; Laplacian matrix;  Critical group;  Invariant factor; Smith
normal form; Spanning tree number.

{\bf 1991 AMS subject classification:}   15A18, 05C50 \\

\noindent{\bf 1. Introduction}

\indent Let $G=(V,E)$ be a finite connected graph without self-loops, but
with multiple edges permitted. Then the Laplacian matrix of $G$ is
the$|V|\times|V|$ matrix defined by
$$L(G)_{uv}
=\left\{\begin{array}{ll}
d(u),&\mbox{if}\quad u=v,\\
\-a_{uv},& \mbox{if}\quad u\neq v,\end{array} \right.$$
where $a_{uv}$ is the number of the edges joining $u$ and $v$, and $d(u)$ is the degree of $u$.\\
\indent Thinking of  $L(G)$ as representing an abelian group homomorphism: $Z^{|V|}\rightarrow Z^{|V|}$, its cokernel has the form
$$Z^{|V|}/\text{im}\, (L(G))\cong Z\oplus Z^{|V|-1}/\mbox{im}\left(\overline{L(G)_{uv}}\right),\eqno(1.1)$$
 where $\overline{L(G)_{uv}}$ is the matrix obtained  from $L(G)$ by striking out  row $u$ and column  $v$,
 and im($\cdot$) refers to the integer span of the columns of the argument. The critical group  $K(G)$ is defined to be
$Z^{|V|-1}/\mbox{im}\left(\overline{L(G)_{uv}}\right)$.
 It is not hard to see that this definition is independent of the choice of $u$ and $v$.
 The critical group $K(G)$ is a finite abelian group, whose order is equal to the absolute value of $\det\overline{L(G)_{uv}}$.
By the well known Kirchhoff's Matrix-Tree Theorem [6, Theorem
13.2.1], the order $|K(G)|$ is equal to the   spanning tree number of $G$.
For the general theory of the critical group, we refer the reader to
Biggs [1, 2], and Godsil [6, Chapter 14].\\
\indent  Recall that an $n\times n$ integral  matrix $P$ is unimodular if $\det P=\pm 1$.
So, the unimodular matrices are precisely those integral matrices with integral inverses, and of course form a
multiplicative group.  Two integral matrices $A$ and $B$ of order $n$ are
equivalent (written by $A \sim B$) if there are unimodular  matrices
$P$ and $Q$ such that $B=PAQ$. Equivalently, $B$ is obtainable from
$A$ by a sequence of elementary row and column operations: (1) the
interchange of two rows or columns, (2) the multiplication of  any
row or column by $-1$, (3) the addition of any integer times of one
row (resp. column) to another row (resp. column). The Smith normal form (Snf) is a diagonal
canonical form for our equivalence relation: every $n\times n$ integral matrix
$A$  is equivalent to a unique diagonal matrix diag$(s_1(A),\cdots,s_n(A))$, where $s_i(A)$ divides $s_{i+1}(A)$ for
$i=1,2,\cdots, n-1$. The $i-$th diagonal entry of the Smith normal form of $A$ is usually called the $i-$th invariant factor of $A$. \\
\indent It is easy to see that $A\sim B$  implies that coker$(A)\cong$ coker$(B)$.
Given any  $n\times n$ unimodular matrices $P$ and $Q$ and any
integral matrix $A$ with $PAQ=$diag$(a_1,\cdots, a_n)$, it is
easy to see that $Z^{|V|}/\mbox{im}(A)\cong (Z/ a_1
Z)\oplus\cdots\oplus(Z/ a_nZ)$.  Assume  the Snf of $\overline{L(G)_{uv}}$
is  diag$(t_1,\cdots,t_{|V|-1})$ (In fact, every such submatrix of $L(G)$ shares the same Snf.), and then it induces an isomorphism
 $$K(G)\cong\left({Z}/ t_1 {Z}\right)\oplus \left({Z}/t_2 {Z}\right)
\oplus\cdots\oplus\left({Z}/ t_{|V|-1} {Z}\right).\eqno(1.2)$$

\indent The nonnegative integers $t_1,t_2,\cdots, t_{|V|-1}$ are also called the invariant factors of $K(G)$, and they can be computed in the following way:
for $1\le i<|V|$, $t_i=\Delta_i/\Delta_{i-1}$
where $\Delta_0=1$ and  $\Delta_i$ is the greatest common divisor of the determinants of the $i\times i$ minors of $\overline{L(G)_{uv}}$. Since $|K(G)|=\kappa$, the spanning tree number of $G$, it follows that $t_1t_2\cdots t_{|V|-1}=\kappa$. So the invariant factors of $K(G)$ can be used to distinguish pairs of non-isomorphic graphs which have the same $\kappa$,
and so there is considerable interest in their properties.
If $G$ is a simple connected graph, the invariant factor $t_1$ of $K(G)$ must be equal to 1,
however, most of them are not easy to be determined.\\

\indent   If $G_1=(V_1, E_1)$ and $G_2=(V_2, E_2)$ are graphs on
disjoint sets of $r$ and $s$ vertices, respectively, their union is
the graph $G_1+G_2=(V_1\cup V_2, E_1\cup E_2)$ and their join
$G_1\vee G_2$ is the graph on $n=r+s$ vertices
obtained from $G_1+G_2$ by inserting new edges joining every vertex
of $G_1$ to every vertex of $G_2$. If we use $G^c$ denote the complement graph of $G$, then  $G_1\vee G_2=(G_1^c+G_2^c)^c$.\\

Compared to the number of the results on the spanning tree number
$\kappa$, there are relatively few results describing the critical
group structure of $K(G)$ in terms of the structure of $G$. There
are also very few interesting infinite family of graphs for which
the group structure has been completely determined (see [3,\, 4,\,
5,\, 7,\, 8,\, 9,\, 11,\, 12] and the references therein). The aim
of this paper is to describe the structure of the critical groups of
two families of graphs $K_m\vee P_n$ and $P_m\vee P_n$}, where
$K_m$ is the complete graph with $m$ vertices, $P_n$ is the path with $n$ vertices.\\


\noindent{\bf\large 2.\, The critical group of $K_m\vee P_n$}\\

\noindent{\bf Lemma 2.1}\, If the graph $G$ has $n$ vertices, then
$$L(K_m\vee G)\sim ((m+1)I_n-L(G^c))\oplus(m+n)I_{m-2}\oplus I_1\oplus 0_1.\eqno(2.1)$$
\noindent{\bf Proof} Note that
$$L(K_m\vee G)=\left(\begin
{array}{cc}{(m+n)I_{m}}-J_m &-J_{m\times{n}}\\-J_{{n}\times
m}&mI_n+L(G)\end{array}\right).$$
\indent Let $$P_1=\left(\begin
{array}{cccccc} 1 & 1 & 1 & \cdots & 1 & 1\\
               1 & 0 & 0 & \cdots & 0 & 0 \\
               -1 & 0 & 1 & \cdots & 0 & 0 \\
                \vdots  & \vdots  & \vdots & \ddots & \vdots &\vdots\\
               -1 & 0 & 0 & \cdots & 1 & 0 \\
               -1 & 0 & 0 & \cdots & 0 & 1
             \end{array}\right),\quad
   Q_1=\left(\begin
{array}{cccccc} 1 & 0 & 0&\cdots & 0 & 0\\
 1 & -1 & -1&\cdots & -1& -1 \\
 1 & 0  & 1 &\cdots &0&0\\
 \vdots&\vdots&\vdots &\ddots &\vdots &\vdots  \\ 1 & 0 & 0&\cdots & 1 & 0 \\ 1 & 0 & 0&\cdots & 0 & 1
             \end{array}\right).$$
Then a direct calculation can
show that $P_1L(K_m\vee G)Q_1=I_1\oplus (m+n)I_{m-2}\oplus((m+n)I_n-L(G^c))\oplus  0_1$.
Note that both the matrices $A$ and $B$ are unimodular, so this Lemma holds.\hfill$\Box$\\

In order to work out the critical group of  graph $K_m\vee P_n$\, $(n\geq4)$,   we only need to work on the Smith normal form of the matrix $(m+n)I_n-L(P_n^c)$.\\
\noindent{\bf Lemma 2.2} $$(m+n)I_n-L(P_n^c)\sim I_{n-2}\oplus\begin{pmatrix}
m+n & b_n \\ 0 &  a_n \\
\end{pmatrix},$$
where
$$\left\{\begin{array}{ll}
a_n&=\dfrac{1}{\sqrt{m^2+4m}}\left(\left(\dfrac{m+2+\sqrt{m^2+4m}}{2}\right)^n
-\left(\dfrac{m+2-\sqrt{m^2+4m}}{2}\right)^n\right),\\
b_n&=e\left(\dfrac{m+2+\sqrt{m^2+4m}}{2}\right)^n
-f\left(\dfrac{m+2-\sqrt{m^2+4m}}{2}\right)^n,
\end{array}\right.$$
and
$$e=\dfrac{(m^2-m-m\sqrt{m^2+4m}-\sqrt{m^2+4m}+2mn)(m+4-\sqrt{m^2+4m})}{4m^2(m+4)},$$
$$f=\dfrac{(m-m^2-m\sqrt{m^2+4m}-\sqrt{m^2+4m}-2mn)(m+4+\sqrt{m^2+4m})}{4m^2(m+4)}.$$

\proof\, Note that
$$(m+n)I_n-L(P_n^c)=\begin
{pmatrix}
{m+2}  &0  &1  &1  &1 & \cdots & 1 &1\\
0     &m+3 &0  &1  &1 & \cdots & 1 &1\\
1    & 0 & m+3 & 0& 1 & \cdots & 1 &1\\
\vdots& \ddots& &\ddots &\ddots & \ddots & &\vdots\\
\vdots& & \ddots& &\ddots&\ddots &\ddots &\vdots\\
1&1&\cdots&1&0&m+3&0&1\\
1&1&\cdots&1&1&0&m+3&0\\
1&1&\cdots&1&1&1&0&m+2
\end{pmatrix}.$$

Let
$P_2=\left(\begin{array}{ccccc}1&0&0&\cdots&0\\-1&1&0&\cdots&0\\
0&-1&1&\cdots&0\\\vdots&\vdots&\ddots&\ddots&\vdots\\
0&0&\cdots&-1&1
\end{array}
\right),$
$Q_2=\left(\begin{array}{ccccc}1&0&\cdots&0&0\\1&1&\cdots&0&0\\
\vdots&\vdots&\ddots&\vdots&\vdots\\1&1&\cdots&1&0\\
1&1&\cdots&1&1
\end{array}
\right),$\\
and  $A_2=P_2((m+n)I_n-L(P_n^c))Q_2$. Then a direct calculation can show
$$A_2=\begin {pmatrix}
{m+n}&{n-2}&{n-2}&{n-3}&n-4&\cdots&1\\
0&m+2&-1&0&0&\cdots&0\\
0&-1&m+2&-1&0&\cdots&0\\
\vdots &\ddots &\ddots&  \ddots  &\ddots &\ddots&\vdots\\
0&\cdots& 0&-1& m+2 & -1&0 \\0&\cdots&\cdots&0&-1&m+2&-1\\
0&\cdots&\cdots&\cdots&0&-1&m+2\\
\end{pmatrix}.$$

For $i=0,\cdots,n-3$, let $M_{i+1}=\begin{pmatrix}
\begin{matrix}
I_{i}&0_{i\times1}\\0_{1\times i}&1\end{matrix}&\begin{matrix} &0_{i\times(n-i-1)}&\\
m+2&-1&0_{1\times(n-i-3)}\end{matrix}\\
0_{(n-i-1)\times (i+1)}&I_{n-i-1} \\
\end{pmatrix}$,\\ and $M_{n-1}=\begin{pmatrix}
\begin{matrix}
I_{n-2}&0_{(n-2)\times1}\\0_{1\times (n-2)}&1\end{matrix}&\begin{matrix} 0_{(n-2)\times1}\\
m+2\end{matrix}\\
0_{1\times (n-1)}&1\\
\end{pmatrix}.$
Let $M=M_1\cdots M_{n-1}$, then $$A_2M=\begin{pmatrix}
\begin{matrix}
m+n\\0\end{matrix}&\begin{matrix}b_2&b_3&b_4&\cdots&b_{n-1}\\
a_2&a_3&a_4&\cdots&a_{n-1}
\end{matrix}&\begin{matrix}b_n\\a_n\end{matrix}\\
0_{(n-2)\times1}&-I_{n-2}&0_{(n-2)\times 1}
\end{pmatrix},$$
where $0_{i\times j}$ is an $i\times j$ zero matrix, and the numbers $a_l$, $b_l$ satisfy the following recurrence relations and initial values
$$\left\{\begin{array}{ll}
a_l=(m+2)a_{l-1}-a_{l-2},& l\ge 3,\\
a_1=1,\quad  a_2=m+2;\\
b_l=(m+2)b_{l-1}-b_{l-2}+(n-l+1),& l\ge 3,\\
b_1=0,\quad  b_2=n-2.
\end{array}\right.\eqno(2.2)$$
Let
$P_3=\begin{pmatrix}
\Large{I_2} & \begin{matrix}
b_2&b_3&\cdots&b_{n-2}&b_{n-1}\\
a_2&a_3&\cdots&a_{n-2}&a_{n-1}\end{matrix}\\
0_{(n-2)\times2}&I_{n-2}
\end{pmatrix},$
Then
$$
P_3A_2M=\begin{pmatrix}
\begin{matrix}m+n\\
0
\end{matrix}&0_{2\times(n-2)}&\begin{matrix}b_n\\
a_n
\end{matrix}\\
0_{(n-2)\times1}&-I_{n-2}&0_{(n-2)\times1}\end{pmatrix}
\sim I_{n-2}\oplus\begin{pmatrix}
m+n & b_n \\ 0 &  a_n \\
\end{pmatrix}.
\eqno(2.3)$$

From (2.2), we can get that
$$a_l=\dfrac{1}{\sqrt{m^2+4m}}\left(\left(\dfrac{m+2+\sqrt{m^2+4m}}{2}\right)^l
-\left(\dfrac{m+2-\sqrt{m^2+4m}}{2}\right)^l\right),$$ and
$$b_l=\left(e\left(\dfrac{m+2+\sqrt{m^2+4m}}{2}\right)^l
-f\left(\dfrac{m+2-\sqrt{m^2+4m}}{2}\right)^l\right)-\dfrac{n-l}{m},$$
where\\
$$e=\dfrac{(m^2-m-(m+1)\sqrt{m^2+4m}+2mn)(m+4-\sqrt{m^2+4m})}{4m^2(m+4)},$$
$$f=\dfrac{(m-m^2-(m+1)\sqrt{m^2+4m}-2mn)(m+4+\sqrt{m^2+4m})}{4m^2(m+4)}.$$
\hfill$\Box$\\

\noindent{\bf Theorem 2.3} (1) The spanning tree number of $K_m\vee P_n$ is
$$\frac{(m+n)^{m-1}}{2^n\sqrt{m^2+4m}}\left(\left(m+2+\sqrt{m^2+4m}\right)^n
-\left(m+2-\sqrt{m^2+4m}\right)^n\right).$$
\indent (2)   The critical group of $K_m\vee P_n$ is
$$Z/(m+n,a_n,b_n) Z\oplus\left(Z/(m+n) Z\right)^{m-2}\oplus Z/\frac{(m+n)a_n}{(m+n,a_n,b_n)}Z,$$
where the parameters $a_n$ and $b_n$ are given in the above Lemma 2.2.\\
\noindent{\bf Proof}  Note that every line sum of the Laplacian matrix of a graph is 0, so we have
 $$L(G)\sim  \mbox{Snf}(\overline{L(G)_{uv}})\oplus 0_1,\quad \mbox{for every }u, v\in V(K_m\vee P_n).\eqno(2.4)$$

It follows from (2.1) and (2.3) that
$$L(K_m\vee P_n)\sim I_{n-1}\oplus\begin{pmatrix}
m+n & b_n \\ 0 &  a_n \\
\end{pmatrix}\oplus(m+n)I_{m-2}\oplus 0_1.\eqno(2.5)$$

Therefore by (2.4) and (2.5), we have
$$\mbox{Snf}(\overline{L(G)_{uv}})\sim I_{n-1}\oplus\begin{pmatrix}
m+n & b_n \\ 0 &  a_n \\
\end{pmatrix}\oplus(m+n)I_{m-2}.\eqno(2.6)$$
It is easy to see that the invariant factors of the matrix $\begin{pmatrix}
m+n & b_n \\ 0 &  a_n \\
\end{pmatrix}\oplus(m+n)I_{m-2}$ are: $(m+n,a_n,b_n)$, $m+n$ (with multiplicity $m-2$), $\dfrac{(m+n)a_n}{(m+n,a_n,b_n)}$, where $(m+n,a_n,b_n)$ stands for the greatest common divisor of $m+n$, $a_n$, and $b_n$.
So this theorem holds.
\hfill$\Box$\\

\noindent{\bf Remark 2.4} It is known that the Laplacian eigenvalues
of $P_n$ are: $0, 2+2\cos{\left(\frac{\pi j}{n}\right)}$ $(1\leq
j\leq n-1)$; and the Laplacian eigenvalues of $K_m$ are: $0, m
$(with multiplicity $m-1$). Then it follows from Theorem 2.1 in [10] that
the  Laplacian eigenvalues of $K_m\vee P_n$ are: $0$, $m+n$ (with
multiplicity $m$), $m+2+2\cos{\left(\frac{j\pi }{n}\right)}$, where
$1\leq j\leq n-1$. Then by the well known Kirchhoff Matrix-Tree
Theorem we know that the spanning tree number of $K_m\vee P_n$ is
$\kappa(K_m\vee
P_n)=(m+n)^{m-1}\prod\limits^{n-1}_{j=1}\left(m+2+2\cos\left(\frac{\pi
j}{n}\right)\right).$ Recall the first part of Theorem 2.3, we have
$$\begin{array}{ll}
&\prod\limits^{n-1}_{j=1}\left(m+2+2\cos\left(\frac{\pi
j}{n}\right)\right)\\
&=\dfrac{1}{2^n\sqrt{m^2+4m}}
\left(\left(m+2+\sqrt{m^2+4m}\right)^n-\left(m+2-\sqrt{m^2+4m}\right)^n\right).
\end{array}\eqno(2.7)$$

\noindent{\bf Example 2.5}\\
\indent If $m=3,n=4$, then $a_4=115,b_4=59$. If $m=4,n=4$, then $a_4=204,b_4=83$. If $m=4,n=5$, then $a_5=1189,b_5=730$.
So it follows from Theorem 2.3 we have the following
$$\begin{array}{lllllllllrrrrr}&&&&&\quad\quad\mbox{Snf}(K_3
\vee P_4)&=&I_4\oplus\mbox{diag}(7,805,0);&&&&&\quad\quad\quad\quad\quad\quad (2.8)\\
&&&&&\quad\quad\mbox{Snf}(K_4\vee P_4)&=&I_4\oplus\mbox{diag}(8,8,1632,0);&&&&&\quad\quad\quad\quad\quad\quad(2.9)\\
&&&&&\quad\quad\mbox{Snf}(K_4 \vee P_5)&=&I_5\oplus\mbox{diag}(9,9,10701,0).&&&&&\quad\quad\quad\quad\quad\quad(2.10)
\end{array}$$
Note that one can use maple to check the results of (2.8), (2.9) and (2.10).\\

\noindent{\bf\Large 3.\, The  critical group of $P_m\vee P_n$}\\

\indent In this section we will work on  the critical group of
$P_m\sq P_n$$(m\geq4,n\geq5)$. Let $L'$ be the submatrix of $L(P_m\sq P_n)$ resulting from the deletion of the last row and the $(m+1)$-th column.
Thus $L'=\left(\begin{array}{cc}
nI_m+L(P_m)&-J_{m\times(n-1)}\\-J_{(n-1)\times{m}}&U\end{array}\right)$,
where $J_{m\times
(n-1)}$ is an $m\times (n-1)$ matrix having all entries equal to 1,
and $U$ is the submatrix obtained from  $mI_n+L(P_n)$ by
deleting its first column and  last row. Now we discuss the Smith
normal form of the matrix $L'$.\\
\indent   Let
$M=\begin{pmatrix}
T_m&0_{m\times(n-1)}\\
0_{(n-1)\times
m}&I_{n-1}
\end{pmatrix}$
and $N=\begin
{pmatrix}T^{-1}_m&0_{m\times(n-1)}\\
0_{(n-1)\times
m}&T^{-1}_{n-1}\end{pmatrix},$
where $T_m=(t_{ij})$ is an $m\times m$ matrix with its entries
$t_{ij}=\left\{\begin{array}{ll}
1,& \mbox{if}\quad i=j, \\
1,&\mbox{if}\quad j=1,  \\
0,     &\mbox{otherwise.}
\end{array}\right.$
Then a direct calculation shows that
$$PNL'M=\begin{pmatrix}
B_{11} & B_{12}\\
B_{21} & B_{22}\end{pmatrix},$$
where  $P$ is an integral matrix with the left multiplication can imply  an  interchange of  row 1 and  $-1$ times  row $m+1$  of $NL'M$, and
$$\begin{array}{ll}
B_{11}&=\begin{pmatrix}
  m     &  1   &  1   &\cdots & \cdots & \cdots &  1\\
  0     &  n+3 &  -1  &\cdots & \cdots & \cdots &  0\\
  0     &  0   &  n+2 & -1    &\cdots  &\cdots  &  0\\
  0     &  1   &  -1  &   n+2 &  -1    &\cdots  &  0\\
 \vdots & \vdots & \vdots & \ddots & \ddots &\ddots &\vdots\\
  0     &  1   &  0   &\cdots &  -1  &  n+2 &  -1  \\
  0     &  1   &  0   &\cdots &\cdots&  -1  &  n+1
\end{pmatrix}_{m\times m},\\
 B_{12}&=\begin{pmatrix}
     \begin{matrix}1  & {0_{1\times(n-2)}}\end{matrix}\\
     0_{(m-1)\times(n-1)} \\
     \end{pmatrix}_{m\times (n-1)},\quad
 B_{21}=\begin{pmatrix}
       \begin{matrix}n  &-1 & 0_{1\times(m-2)}\end{matrix}\\
       0_{(n-2)\times m}\\
  \end{pmatrix}_{(n-1)\times m}, \\
  B_{22}&=\begin{pmatrix}
    -1   &   -1    &   -1   &\cdots &  \cdots &   -1  &   -1 \\
  m+3   &  -1    &   0     &\cdots  &  \cdots &   0   &   0  \\
   0     &  m+2   &  -1     &\cdots & \cdots &   0   &   0  \\
   1     &  -1    &  m+2    & -1    & \cdots &   0   &   0  \\
 \vdots  & \vdots & \ddots &\ddots &\ddots &\vdots &\vdots \\
   1     &   0    &\cdots& -1   & m+2  &  -1  &   0  \\
   1     &   0    &\cdots&\cdots&  -1  &  m+2 &  -1
 \end{pmatrix}_{(n-1)\times (n-1)}. \end{array} $$

The sequences $p_i$, $q_i$, $c_i$, $d_i$ will be used in the
following  Lemma 3.1, where

$$\left\{
\begin{array}{lll}
  p_k=(n+2)p_{k-1}+q_{k-1},&k\geq1\\
  q_k=1-p_{k-1},\\
  p_0=n+3,\, q_0=0.
\end{array}
\right.\eqno (3.1)
$$
and $$\left\{
\begin{array}{lll}
  c_k=(m+2)c_{k-1}+d_{k-1}, & k\geq1\\
  d_k=1-c_{k-1},\\
  c_0=m+3,\, d_0=0.
\end{array}
\right.\eqno (3.2)
$$

\noindent{\bf Lemma 3.1}
$$L'\sim F=\left(\begin {array}{ccc}
I_{(m+n-3)}& 0_{(m+n-3)\times2}\\
\Large{0_{2\times(m+n-3)}}&\begin{matrix}p'_{m-2}&0\\
\alpha\beta-1&m\beta+n\\
\end{matrix}\end{array}
\right),$$
where
$$p'_{m-2}=\dfrac{1}{\sqrt{n^2+4n}}\left(\left(\dfrac{n+2+\sqrt{n^2+4n}}{2}\right)^m
-\left(\dfrac{n+2-\sqrt{n^2+4n}}{2}\right)^m\right),$$
$$\alpha=\dfrac{1}{n}{p'_{m-2}}-\dfrac{m}{n},$$
$$\beta=\dfrac{1}{m\sqrt{m^2+4m}}\left(\left(\dfrac{m+2+\sqrt{m^2+4m}}{2}\right)^n
-\left(\dfrac{m+2-\sqrt{m^2+4m}}{2}\right)^n\right)-\frac{n}{m}.$$

\proof\,
Let $$
M_1=\begin{pmatrix}
m  &  1  &  1  &  1  &\cdots&\cdots&  1  \\
0  & p_0 & -1  &  0  &\cdots&\cdots&  0  \\
0  & q_0 & n+2 & -1  &\cdots&\cdots&  0  \\
0  & 1   & -1  & n+2 & -1   &\cdots&  0  \\
\vdots&\vdots&\vdots&\ddots&\ddots&\ddots&\vdots\\
0  &  1  &  0  &\cdots& -1 & n+2 &  -1 \\
0  &  1  &  0  &\cdots&\cdots  & -1  & n+1
\end{pmatrix}_{m\times m}
$$
and
$$
M_2=\left(
\begin{array}{cccccccccccc}
 -1   &   -1    &   -1   &\cdots &  \cdots &   -1  &   -1 \\
c_0   &  -1    &   0     &\cdots  &  \cdots &   0   &   0  \\
d_0     &  m+2   &  -1     &\cdots & \cdots &   0   &   0  \\
   1     &  -1    &  m+2    & -1    & \cdots &   0   &   0  \\
 \vdots  & \vdots & \ddots &\ddots &\ddots &\vdots &\vdots \\
   1     &   0    &\cdots& -1   & m+2  &  -1  &   0  \\
   1     &   0    &\cdots&\cdots&  -1  &  m+2 &  -1
\end{array}
\right)_{(n-1)\times (n-1)},
$$
then we can rewrite $L'$ as~$L'=\begin{pmatrix}
          M_1 & B_{12} \\
          B_{21} & M_2
        \end{pmatrix}.$

Let $B=\left(
     \begin{array}{cc}
       B_1 & 0 \\
       0 & B_2 \\
     \end{array}
   \right)
$
where
$$B_1=\begin{pmatrix}
      1 & 0 & 0 & 0 & \cdots & 0 \\
      0 & 1 & 0 & 0 & \cdots & 0 \\
      0 & p_0 & 1 & 0 & \cdots & 0 \\
      0 & p_1 & 0 & 1 & \cdots & 0 \\
      \vdots & \vdots & \vdots & \vdots & \ddots & \vdots \\
      0 & p_{m-3} & 0 & 0 & \cdots & 1 \\
    \end{pmatrix}_{m\times m},\quad B_2=
    \begin{pmatrix}
       1 & 0 & 0 & \cdots & 0 \\
       c_0 & 1 & 0 & \cdots & 0 \\
       c_1 & 0 & 1 & \cdots & 0 \\
      \vdots & \vdots & \vdots & \ddots & \vdots \\
       c_{n-3} & 0 & 0 & \cdots & 1 \\
    \end{pmatrix}_{{(n-1)}\times {(n-1)}}.$$

Then it is easy to check that
$$L'B=\begin{pmatrix}
          M'_1 & B_{12} \\
          B_{21} & M'_2 \\
        \end{pmatrix},$$
where
$$
M'_1=\begin{pmatrix}
m  &  \alpha  &  1  &  1  &  \cdots  & \cdots  &  1  \\
0  & 0 & -1  &  0  &  \cdots  &  \cdots  &  0  \\
0  &  0  & n+2 & -1 &  \cdots  & \cdots  &  0  \\
0  &  0  & -1 & n+2 & -1  & \cdots  &  0  \\
\vdots&\vdots&\vdots&\ddots&\ddots&\ddots&\vdots\\
0  &  0  & 0 & \cdots &  -1  & n+2 &  -1 \\
0  &  p'_{m-2}  & 0 & \cdots &  \cdots  & -1  & n+1 \\
\end{pmatrix}_{m\times m},
$$
and $p'_{m-2}=(n+1)p_{m-3}-p_{m-4}+1$,
$\alpha=1+\sum\limits_{k=0}^{m-3}p_k.$

From (3.1), we get $$p_k=(x+y)\left(\dfrac{n+2+\sqrt{n^2+4n}}{2}\right)^k
+(x-y)\left(\dfrac{n+2-\sqrt{n^2+4n}}{2}\right)^k-\dfrac{1}{n},$$ where $x=\dfrac{n^2+3n+1}{2n}$, $y=\dfrac{n^2+5n+5}{2\sqrt{n^2+4n}}.$ And now we can easily get the expression of $p'_{m-2}$ and $\alpha$ by a direct calculation.
$$
M'_2=
\begin{pmatrix}
-\beta &  -1 &  -1 &  -1 &  \cdots  & -1  & -1  \\
0&  -1 &  0  &  0  &  \cdots  &  0  &  0  \\
0  & m+2 & -1  &  0  &  \cdots  &  0  &  0  \\
0  & -1  & m+2 & -1  &  \cdots  &  0  &  0  \\
\vdots&\vdots&\ddots&\ddots&\ddots&\vdots&\vdots\\
0  &  0  &\cdots & -1 &  m+2 & -1  &  0  \\
0  &  0  &\cdots &\cdots &  -1  & m+2 & -1  \\
\end{pmatrix}_{{(n-1)}\times {(n-1)}},
$$
and $\beta=1+\sum\limits_{k=0}^{n-3}c_k.$

From (3.2), we can obtain that $$c_k=(u+v)\left(\dfrac{m+2+\sqrt{m^2+4m}}{2}\right)^k
+(u-v)\left(\dfrac{m+2-\sqrt{m^2+4m}}{2}\right)^k-\dfrac{1}{m},$$ where $u=\dfrac{m^2+3m+1}{2m}$, $v=\dfrac{m^2+5m+5}{2\sqrt{m^2+4m}}$. Thus we get $\beta$ by a direct calculation.

Now we deal with the matrix $L'B$. In the following, we will use $r_i$ to denote the $i-$th row of matrix $L'B$.

For $i=2,\cdots,m-2$, we first add $(n+2)r_i$ to $r_{i+1}$ and add $-r_i$ to $r_{i+2}$; then add $r_i$ to $r_1$, and
add $(n+1)r_{m-1}$ to $r_m$. After that we have $$M'_1\sim M''_1=
\begin{pmatrix}
\begin{matrix}m  &  \alpha  &  0_{1\times(m-2)}\end{matrix}\\
\begin{matrix}0_{(m-2)\times2} &-I_{(m-2)} \end{matrix}\\
\begin{matrix}0  & p'_{m-2} &  0_{1\times(m-2)}\end{matrix}\\
\end{pmatrix}_{m\times m}.$$

For $i=m+2,\cdots,m+n-2$, we first add $(m+2)r_i$ to $r_{i+1}$ and add $-r_i$ to $r_{i+2}$; then add $-r_i$ to $r_{m+1}$, and add $-r_{m+n-1}$ to $r_{m+1}$. After that we have $$M'_2\sim M''_2=
\begin{pmatrix}
-\beta& 0_{1\times(n-2)}\\
0_{(n-2)\times1}&-I_{(n-2)}
\end{pmatrix}_{(n-1)\times (n-1)}.$$
Note that the matrices $B_{12}$ and $B_{21}$ are not influenced in the operations on $L'B$.

So
$$\begin{array}{ll}L'B & \sim\begin{pmatrix}
          M''_1 & B_{12} \\
          B_{21} & M''_2 \\
        \end{pmatrix}
\sim I_{m+n-4}\oplus\begin{pmatrix}
     m & \alpha & 1 \\
     0 & p'_{m-2} & 0 \\
     n & -1 & -\beta
   \end{pmatrix}\\
&\sim I_{m+n-4}\oplus\begin{pmatrix}
     m & \alpha & 1 \\
     0 & p'_{m-2} & 0 \\
     m\beta+n & \alpha\beta-1 & 0
   \end{pmatrix}\\
&\sim I_{m+n-3}\oplus
   \begin{pmatrix}
 p'_{m-2}  &  0      \\
 \alpha\beta-1     &  m\beta+n
   \end{pmatrix}.\end{array}$$\hfill$\Box$\\

From lemma 3.2, we immediately obtain the following theorem.

\noindent{\bf Theorem 3.2}
  The critical group of the graph $P_m\sq P_n$ is $Z/tZ\oplus Z/sZ$,
where $t=(m\beta+n,\alpha\beta-1,p'_{m-2}), s=\dfrac{(m\beta+n)
p'_{m-2}}{(m\beta+n,\alpha\beta-1,p'_{m-2})}$, where $\alpha$ and $\beta$ are
defined as above. The spanning tree number of  $P_m\sq P_n$ is
$(m\beta+n)p'_{m-2}$.\\

\noindent{\bf Remark 3.3} It is known that the Laplacian eigenvalues
of $P_m$ are: $0, 2+2\cos{\left(\frac{\pi j}{m}\right)}$, $(1\leq j\leq n-1)$.
Then it follows from Theorem 2.1
in [10] that the Laplacian eigenvalues of $P_m\vee P_n$ are: $0$,
$m+n$, $n+2+2\cos{\left(\frac{i\pi}{m}\right)}(1\leq i\leq m-1)$,
$m+2+2\cos{\left(\frac{j\pi }{n}\right)}(1\leq j\leq n-1)$. Then by
the well known Kirchhoff Matrix-Tree Theorem we know that the
spanning tree number of $P_m\vee P_n$ is $\kappa(P_m\vee
P_n)=\prod\limits^{m-1}_{i=1}\left(n+2+2\cos{\left(\frac{i\pi}{m}\right)}\right)\prod\limits^{n-1}_{j=1}\left(m+2+2\cos{\left(\frac{j\pi
}{n}\right)}\right).$ Recall Theorem 3.2, we have
$$\begin{array}{ll}
&\prod\limits^{m-1}_{i=1}\left(n+2+2\cos{\left(\frac{i\pi}{m}\right)}\right)\prod\limits^{n-1}_{j=1}\left(m+2+2\cos{\left(\frac{j\pi
}{n}\right)}\right)\\
&=\frac{\left(\left(n+2+\sqrt{n^2+4n}\right)^m
-\left(n+2-\sqrt{n^2+4n}\right)^m\right)
\left(\left(m+2+\sqrt{m^2+4m}\right)^n-\left(m+2-\sqrt{m^2+4m}\right)^n\right)}{2^{m+n}\sqrt{(n^2+4n)(m^2+4m)}}.
\end{array}\eqno(3.3)$$

\noindent{\bf Example 3.4}\\
\indent If $m=4,n=5$, then $p'_3=329, \alpha=65, \beta=296$. If $m=4,n=6$,
then $p'_4=496,\alpha=82, \beta=1731$. If $m=5,n=5$, then
$p'_3=2255, \alpha=450, \beta=450$. So it follows from Theorem 3.2 that
we have the following
$$\begin{array}{lllllllllrrrrr}&&&&&\quad\quad\mbox{Snf}(P_4 \vee P_5)&=&I_7\oplus\mbox{diag}(391181,0);&&&&&&\quad\quad\quad\quad\quad\quad (3.4)\\
&&&&&\quad\quad\mbox{Snf}(P_4 \vee P_6)&=&I_8\oplus\mbox{diag}(3437280,0);&&&&&&\quad\quad\quad\quad\quad\quad(3.5)\\
&&&&&\quad\quad\mbox{Snf}(P_5 \vee P_5)&=&I_7\oplus\mbox{diag}(451,11275,0).&&&&&&\quad\quad\quad\quad\quad\quad(3.6)
\end{array}$$
Here we also note that one can use Maple to check the results of (3.4), (3.5) and (3.6).\\

\end{document}